\magnification=1200

\centerline{\bf Decoupling Inequalities for the Tail Probabilities of
Multivariate U-statistics}

\centerline{\bf by}

\centerline{Victor H. de la Pe\~na$^{1}$ and S. J. Montgomery-Smith $^{2}$}

\centerline{Columbia University and University of Missouri, Columbia}
\vskip 0.2in
\centerline{\bf Abstract}
\vskip 0.2in

In this paper we present a decoupling inequality that shows that
multivariate U-statistics can be studied as sums of (conditionally)
independent random variables. This result has important implications
in several areas of probability and statistics including
the study random graphs and multiple stochastic integration. 
More precisely, we get the following result:

\noindent {\bf Theorem 1. }  Let $\{X_j\}$  be a sequence of independent
random variables in a measurable space $({\cal S},S)$, and let $\{ X_i^{(j)}\}, ~j=1,...,k$
be $k$ independent copies of $\{X_i\}$. Let $f_{i_1i_2...i_k}$ be families of 
functions of $ k $ variables taking $(S \times ... \times S)$ into a Banach
space $(B,||\cdot||)$. Then, for all $n\ge k \ge 2,~t>0,$ there exist numerical
constants $C_k$ depending on $k$ only so that,

$$P(||\sum_{1\le i_1 \ne i_2 \ne ... \ne i_k \le n} f_{i_1...i_k}
(X_{i_1}^{(1 )},X_{i_2}^{(1)},...,X_{i_k}^{(1)}) || \ge t) 
$$

$$ \le C_k
P( C_k||\sum_{1\le i_1 \ne i_2 \ne ... \ne i_k \le n} f_{i_1...i_k}
(X_{i_1}^{(1 )},X_{i_2}^{(2)},...,X_{i_k}^{(k)}) || \ge t).
$$

The reverse bound holds if in addition, the following symmetry condition 
holds almost surely
$$f_{i_1i_2...i_k}(X_{i_1},X_{i_2},...,X_{i_k}) =
f_{i_{ \pi (1)}i_{\pi (2)}...i_{\pi (k)}}(X_{ i_{\pi (1)}},X_{i_{\pi (2)}},...,
X_{ i_{\pi (k)}}),$$ 
\noindent for all permutations $\pi $ of $(1,...,k)$.

\centerline{\bf 1. Introduction}

In this paper we provide the multivariate extension of the  tail probability
decoupling
inequality for generalized U-statistics of order two and quadratic forms
presented in de la Pe\~na and Montgomery-Smith (1993). This type of 
inequality permits the transfer of some results
for sums of independent random variables to the case of U-statistics.
Our work  builds mainly on recent work of 
Kwapien and Woyczynski (1992) as well as on results for U-statistics from  
Gin\'e  and Zinn (1992) 
and papers dealing with inequalities for multilinear forms
of symmetric and
hypercontractive random variables in 
de la Pe\~na, Montgomery-Smith and Szulga (1992), and de la Pe\~na (1992).
It is to be remarked that the decoupling inequalities for multilinear forms
introduced in McConnell and Taqqu (1986) provided us with our first 
exposure to this decoupling problem. For a more expanded list of
references on the subject
see, for example, Kwapien and Woyczynski (1992).

\vskip 0.2in
\noindent $^{1,2}$ Supported in part by NSF grants. 
\par
\noindent $^2$ Supported by the University of Missouri Research Board.



\eject

\centerline{\bf 2. Main Result}

\noindent {\bf Theorem 1.}  Let $\{X_i\}$  be a sequence of independent
random variables in a measurable space $({\cal S},S)$, and let $\{ X_i^{(j)}\}, ~j=1,...,k$
be $k$ independent copies of $\{X_i\}$. Let $f_{i_1i_2...i_k}$ be families of 
functions of $ k $ variables taking $(S \times ... \times S)$ into a Banach
space $(B,||\cdot||)$. Then, for all $n\ge k\ge 2,~t>0,$ there exist numerical
constants $C_k,~\tilde C_k$ depending on $k$ only so that,
$$P(||\sum_{1\le i_1 \ne i_2 \ne ... \ne i_k \le n} f_{i_1...i_k}
(X_{i_1}^{(1 )},X_{i_2}^{(1)},...,X_{i_k}^{(1)}) || \ge t) 
$$

$$ \le C_k
P( C_k||\sum_{1\le i_1 \ne i_2 \ne ... \ne i_k \le n} f_{i_1...i_k}
(X_{i_1}^{(1 )},X_{i_2}^{(2)},...,X_{i_k}^{(k)}) || \ge t).
$$

If in addition, the following symmetry condition holds almost surely
$$f_{i_1i_2...i_k}(X_{i_1},X_{i_2},...,X_{i_k}) =$$
$$f_{i_{\pi (1)}i_{\pi (2)}...i_{\pi (k)}}(X_{i_{\pi (1)}},X_{i_{\pi (2)}}
,...,
X_{i_{\pi (k)}})$$
for all permutations $\pi$ of $(1,...,k)$, then 

$$P(||\sum_{1\le i_1 \ne i_2 \ne ... \ne i_k \le n} f_{i_1...i_k}
(X_{i_1}^{(1)},X_{i_2}^{(2)},...,X_{i_k}^{(k)}) || \ge t) 
$$

$$ \le \tilde C_k
P(\tilde C_k||\sum_{1\le i_1 \ne i_2 \ne ... \ne i_k \le n} f_{i_1...i_k}
(X_{i_1}^{(1 )},X_{i_2}^{(1)},...,X_{i_k}^{(1)}) || \ge t).
$$

\noindent{\bf Note:} In this paper we use the notation
$\{ i_1 \ne i_2 \ne ...\ne i_k \}$ to denote that all of
$i_1,...,i_k$ are different.

\centerline{\bf 3. Preliminary Results}

Throughout this paper we will be using two results found in
earlier work. The first one comes from de la Pe\~na and 
Montgomery-Smith (1993).
For completeness we reproduce the proof here.

\noindent{\bf Lemma 1.}
Let $X,~Y$ be two i.i.d. random variables. Then
$$P(||X||\ge t) \le 3P(||X+Y||\ge{2t\over 3}).\leqno (1)$$
\noindent{\bf Proof:} Let $X, Y, Z$ be i.i.d. random
variables. Then
$$\eqalign{P(||X ||& \ge t)  \cr
& = P(|| (X+Y) + (X+Z) - (Y+Z) || \ge 2t) \cr
& \le
P(||X+Y||\ge 2t/3) +  P(||X+Z||\ge 2t/3) +  P(||Y+Z||\ge 2t/3) \cr
& = 
3 P(||X+Y||\ge 2t/3) .\cr}$$

The second result comes from Kwapien and Woyczynski (1992) and can also be
found in de la Pe\~na and Montgomery-Smith (1993).

\noindent{\bf Proposition 1.} Let $Y$ be any mean zero random variable 
with values in  a Banach space $(B,||\cdot||)$. Then, for all $a\epsilon B$,
$$P(||a+Y||\ge||a||) \ge {\kappa \over 4},\leqno (2)$$
where $\kappa =
\inf_{x'\epsilon B'}{(E|x'(Y)|)^2 \over E(x'(Y))^2}.$
(Here $B'$\ denotes the family of linear functionals on $B$.)

\noindent{\bf Proof:}  Note first that if $\xi$\ is a random variable
for which $E\xi = 0$, then
$ P (\xi \ge 0 ) \ge {1\over 4} {(E|\xi |)^2 \over E(\xi ^2) } $. 
From this, we deduce that
$ P(x'(Y) \ge 0) \ge {1\over 4}{(E|X'(Y)|)^{2}\over E(X'(Y))^{2}}$ 
The result then follows, because if
$x' \epsilon B'$ 
is such that
$||x'||=1$ and $x'(a) =||a||$, then $\{||a+Y||\ge||a||\}$ contains
$\{x'(a+Y) \ge x'(a)\} = \{x'(Y) \ge 0\}$. 

\noindent{\bf Lemma 2.} Let $x, a_{i_1}, a_{i_1i_2},...,a_{i_1i_2...i_k}$ 
belong
to a Banach space $(B,||\cdot ||)$.
Let $\{\epsilon_i\}$ be a sequence of 
symmetric Bernoulli random variables. Then,

$$P(||x+\sum_{r=1}^k \sum_{1\le i_i \ne i_2 \ne ... \ne i_r \le n} 
a_{i_1...i_r} 
\epsilon_{i_1} ...\epsilon_{i_r}|| \ge
||x||) \ge c_k ^{-1},$$
for a universal constant $1<c_k<\infty$ depending on $k$ only.

\noindent{\bf Proof:} Suppose that 
$x, a_{i_1}, a_{i_1i_2},...,a_{i_1i_2...i_k}$ are in $R$, then
since the $\epsilon$'s are hypercontractive, by equation (1.4) of
Kwapien and Szulga (1991) and the easy argument of the proof of Lemma 3
in de la Pe\~na and Montgomery-Smith (1993), for some $\sigma > 0$, we get

\vskip 0.1in
\noindent $(E|\sum_{r=1}^k \sum_{1\le i_1 \ne ... \ne  i_r \le n} 
a_{i_1...i_r}
\epsilon_{i_1}...\epsilon_{i_r}|^4)^{1\over 4}  $

\noindent $=(E|\sum_{r=1}^k \sum_{1\le i_1 < ...  < i_r \le n}
b_{i_1...i_r}
\epsilon_{i_1}...\epsilon_{i_r}|^4)^{1\over 4}  $

\noindent $ \le \sigma^{-k}
(E|\sum_{r=1}^k \sum_{1\le i_1 < ... < i_k \le n} b_{i_1...i_r}
\epsilon_{i_1}...\epsilon_{i_r}|^2)^{1\over 2}  $

\noindent $ = \sigma^{-k}
(E|\sum_{r=1}^k \sum_{1\le i_1 \ne  ... \ne i_k \le n} a_{i_1...i_r}
\epsilon_{i_1}...\epsilon_{i_r}|^2)^{1\over 2},  $
\vskip 0.1in
\noindent where $b_{i_1...i_r} = \sum_{\pi\in S_r} a_{i_{\pi(1)}...i_{\pi(r)}}$,
and $S_r$\ denotes the set of all permutations of $\{1,...,r\}$.

Next, observe that $||\xi||_4 \le \sigma^{-2} ||\xi ||_2 $ implies that
$||\xi ||_2 \le \sigma^{-4} ||\xi ||_1$. Take $x' \epsilon B'$ so that
$||x'||=1 $ and $x'(x) = ||x||$, then

\vskip 0.1in
\noindent $P(||x+\sum_{r=1}^k \sum_{1\le i_i \ne i_2 \ne ... \ne i_r \le n} 
a_{i_1...i_r} 
\epsilon_{i_1} ...\epsilon_{i_r}|| \ge
||x||) $

\noindent $\ge 
P(x'(x) + \sum_{r=1}^k \sum_{1\le i_i \ne i_2 \ne ... \ne i_r \le n} 
x'(a_{i_1...i_r} )
\epsilon_{i_1} ...\epsilon_{i_r} \ge
x'(x)) $

\noindent $ =
P(\sum_{r=1}^k \sum_{1\le i_i \ne i_2 \ne ... \ne i_r \le n} 
x'(a_{i_1...i_r} )
\epsilon_{i_1} ...\epsilon_{i_r} \ge
0) \ge c_k^{-1}$
\vskip 0.1in
\noindent {\bf Note:}
Throughout this paper we will use $c_k$ and $C_k$ to denote 
numerical constants that depend on
$k$ only and may change from application to application. 

\centerline{\bf 4. Proof of the Upper Bound:}

Our proof of this result is obtained
by applying the argument used in the proof of 
the upper bound in the bivariate case plus an inductive argument.
Let $\{\sigma_i\}$ be a sequence of independent
symmetric Bernoulli random variables, $P(\sigma_i=1) = {1\over 2} $ and 
$P(\sigma_i=-1) ={ 1\over 2}$.
Consider random variables $(Z_i^{(1)},Z_i^{(2)})$
such that $(Z_i^{(1)},Z_i^{(2)}) = (X_i^{(1)},X_i^{(2)})$ if 
$\sigma_i =1$ and
$(Z_i^{(1)},Z_i^{(2)}) = (X_i^{(2)},X_i^{(1)})$ if $\sigma_i =-1$.
Then $(1+\sigma_i)$ and $(1-\sigma_i)$ are either $0$ or $2$ and these random
variables can be used to transform the problem from one involving $X$'s to 
one involving $Z$'s.  Let us first illustrate the argument in the case that
$k=3$.

\noindent $2^3f_{i_1 i_2 i_3}(Z_{i_1}^{(1)}, Z_{i_2}^{(1)}, Z_{i_3}^{(2)}) = $

$$\{(1+\sigma_{i_1})(1+\sigma_{i_2})(1+\sigma_{i_3})
f_{i_1 i_2 i_3}(X_{i_1}^{(1)}, X_{i_2}^{(1)}, X_{i_3}^{(2)})  $$ 

$$+(1+\sigma_{i_1})(1+\sigma_{i_2})(1-\sigma_{i_3})
f_{i_1 i_2 i_3}(X_{i_1}^{(1)}, X_{i_2}^{(1)}, X_{i_3}^{(1)})  $$

$$ +(1+\sigma_{i_1})(1-\sigma_{i_2})(1+\sigma_{i_3})
f_{i_1 i_2 i_3}(X_{i_1}^{(1)}, X_{i_2}^{(2)}, X_{i_3}^{(2)})  $$

$$+(1-\sigma_{i_1})(1+\sigma_{i_2})(1+\sigma_{i_3})
f_{i_1 i_2 i_3}(X_{i_1}^{(2)}, X_{i_2}^{(1)}, X_{i_3}^{(2)}) \leqno (3) $$

$$+(1+\sigma_{i_1})(1-\sigma_{i_2})(1-\sigma_{i_3})
f_{i_1 i_2 i_3}(X_{i_1}^{(1)}, X_{i_2}^{(2)}, X_{i_3}^{(1)})  $$

$$+(1-\sigma_{i_1})(1+\sigma_{i_2})(1-\sigma_{i_3})
f_{i_1 i_2 i_3}(X_{i_1}^{(2)}, X_{i_2}^{(1)}, X_{i_3}^{(1)})  $$ 

$$+(1-\sigma_{i_1})(1-\sigma_{i_2})(1+\sigma_{i_3})
f_{i_1 i_2 i_3}(X_{i_1}^{(2)}, X_{i_2}^{(2)}, X_{i_3}^{(2)})  $$

$$ +(1-\sigma_{i_1})(1-\sigma_{i_2})(1-\sigma_{i_3})
f_{i_1 i_2 i_3}(X_{i_1}^{(2)}, X_{i_2}^{(2)}, X_{i_3}^{(1)})\},   $$

\noindent  
where the sign ``$+$'' is chosen if the superscript
of $X_i$\ agrees with that of $Z_i$, and ``$-$'' otherwise.
Next, set $  T_{n,3} = \displaystyle{\sum_{1\le i_1 \ne i_2 \ne i_3 \le n}}$
$$ \{f_{i_1i_2 i_3}(X_{i_1}^{(1)}, X_{i_2}^{(1)}, X_{i_3}^{(2)})  
+f_{i_1 i_2 i_3}(X_{i_1}^{(1)}, X_{i_2}^{(1)}, X_{i_3}^{(1)})  $$
$$+f_{i_1 i_2 i_3}(X_{i_1}^{(1)}, X_{i_2}^{(2)}, X_{i_3}^{(2)}) 
+f_{i_1 i_2 i_3}(X_{i_1}^{(2)}, X_{i_2}^{(1)}, X_{i_3}^{(2)})  $$
$$
+f_{i_1 i_2 i_3}(X_{i_1}^{(1)}, X_{i_2}^{(2)}, X_{i_3}^{(1)}) + 
f_{i_1 i_2 i_3}(X_{i_1}^{(2)}, X_{i_2}^{(1)}, X_{i_3}^{(1)}) $$
$$+f_{i_1 i_2 i_3}(X_{i_1}^{(2)}, X_{i_2}^{(2)}, X_{i_3}^{(2)}) + 
f_{i_1 i_2 i_3}(X_{i_1}^{(2)}, X_{i_2}^{(2)}, X_{i_3}^{(1)})\}. $$ 

Letting ${ \cal G}_2 = \sigma (X_i^{(1)},X_i^{(2)}, i=1,...,n)$ we get
$$T_{n,3}= 
2^3 \sum_{1\le i_1 \ne i_2 \ne i_3 \le n} 
E(f_{i_1 i_2 i_3}(Z_{i_1}^{(1)}, Z_{i_2}^{(1)}, Z_{i_3}^{(2)})| {\cal G}_2).  $$

More generally, for any $1\le l_1,...,l_k \le 2$, one can obtain the 
expansion 

\noindent $2^k f_{i_1...i_k}(Z_{i_1}^{(l_1)},...,Z_{i_k}^{(l_k)}) $
$$= \sum_{1\le j_1,...,j_k \le 2} (1{\pm}\sigma_{i_1})...(1{\pm}
\sigma_{i_k})f_{i_1...i_k}(X_{i_1}^{(j_1)},...,X_{i_k}^{(j_k)}). \leqno
(4)$$

The appropriate extension of $T_{n,3}$ is 
$$T_{n,k} = \sum_{1\le i_1 \ne ... \ne i_k \le n} 
\sum_{1\le j_1,...,j_k \le 2} f_{i_1...i_k}
(X_{i_1}^{(j_1)},...,X_{i_k}^{(j_k)}).$$
Again, 
$$T_{n,k} = 2^k\sum_{1\le i_1 \ne ... \ne i_k \le n} E(f_{i_1...i_k}
(Z_{i_1}^{(l_1)},...,Z_{i_k}^{(l_k)})|{\cal G}_2).$$

From Lemma 1 
we get, 
$$P(||\sum_{1\le i_1 \ne i_2 \ne ... \ne i_k \le n} f_{i_1...i_k}
(X_{i_1}^{(1 )},...,X_{i_k}^{(1)}) || \ge t) \le 
$$
$$ 3P(3||\sum_{1\le i_1 \ne ... \ne i_k \le n} \{f_{i_1...i_k}
(X_{i_1}^{(1 )},...,X_{i_k}^{(1)}) +
f_{i_1...i_k}
(X_{i_1}^{(2 )},...,X_{i_k}^{(2)})\}
 || \ge 2 t) = $$
$$ 3P(3||T_{n,k} +\sum_{1\le i_1 \ne ... \ne i_k \le n} f_{i_1...i_k}
(X_{i_1}^{(1 )},...,X_{i_k}^{(1)}) +
f_{i_1...i_k}
(X_{i_1}^{(2 )},...,X_{i_k}^{(2)}) - T_{n,k} || \ge 2t) \le $$
\noindent $\{3P(3||T_{n,k}|| \ge t) $
$$+ 3P(3||\sum_{1\le i_1 \ne ... \ne i_k\le n} 
\sum_{1\le j_1,...j_k \le 2,  \hbox {  not all j's equal} }
f_{i_1 ... i_k}(X_{i_1}^{(j_1)},...,X_{i_k}^{(j_k)})|| \ge t)\} $$
\noindent $\le \{3P(3||T_{n,k}|| \ge t) $
$$+ 
\sum_{1\le j_1,...,j_k \le 2,  \hbox { not all j's equal} }
C_kP(C_k||\sum_{1\le i_1\ne ...\ne i_k\le n} 
f_{i_1 ... i_k}(X_{i_1}^{(j_1)},...,X_{i_k}^{(j_k)})||\ge t)\}. 
\leqno (5)$$
(Recall that $C_k,c_k$ are numerical constants that depend on $k$ only
and may change from application to application.)

Observe also that using (4) and the fact that the $\sigma$'s are
independent from the $X$'s, 
Lemma 2 with $x=T_{n,k}$
gives
for any fixed
$1\le l_1,...,l_k \le 2$,
$$P(2^k||\sum_{1\le i_1 \ne ... \ne i_k\le n}
f_{i_1 ... i_k}(Z_{i_1}^{(l_1)},...
,Z_{i_k}^{(l_k)})|| \ge ||T_{n,k}|| ~|{\cal G}_2)\ge c_k^{-1},\leqno (6)$$

Integrating over $\{||T_{n,k}||\ge t\}$ and using the fact that
$\{(X_i^{(1)},X_i^{(2)}): i=1,...,n \}$ has the 
same joint distribution as 
$\{(Z_i^{(1)},Z_i^{(2)}): i=1,...,n \}$ we obtain that 
$$P(2^k||\sum_{1\le i_1 \ne ... \ne i_k \le n}
f_{i_1 ... i_k}(X_{i_1}^{(l_1)},...
,X_{i_k}^{(l_k)})|| \ge t ). \leqno (7)$$
$$=P(2^k||\sum_{1\le i_1 \ne ... \ne i_k\le n}f_{i_1 ... i_k}(Z_{i_1}^{(l_1)},
...
,Z_{i_k}^{(l_k)})|| \ge t) \ge  c_k^{-1}P(||T_{n,k}||\ge t)$$

It is obvious that the upper bound decoupling 
inequality holds for the case of U-statistics of order $1$. 
Assume that it holds for U-statistics of orders
$2,...,k-1$.
Putting (5) and (7) together with $1\le l_1,...,l_k \le 2$, not all
$l$'s equal  
we get,

\noindent $P(||\sum_{1\le i_1 \ne ... \ne i_k \le n} f_{i_1...i_k}
(X_{i_1}^{(1 )},...,X_{i_k}^{(1)}) || \ge t) $
\noindent $\le \{3P(3||T_{n,k}|| \ge t)$
$$+ 
\sum_{1\le j_1,...,j_k \le 2, \hbox { not all j's equal} }
C_kP(C_k||\sum_{1\le i_1\ne ... \ne i_k\le n} 
f_{i_1 ...i_k}(X_{i_1}^{(j_1)}, ...,X_{i_k }^{(j_k)})||\ge t)
\}$$
$$\le \sum_{1\le j_1,...,j_k \le 2, \hbox{ not all j's equal}}
C_kP(C_k||\sum_{1\le i_1 \ne ...  \ne i_k \le n} f_{i_1...i_k}
(X_{i_1}^{(j_1 )},...,X_{i_k}^{(j_k)}) || \ge t) $$

\noindent $ ~~~~~\le C_kP(C_k||\sum_{1\le i_1 \ne ... \ne i_k \le n} 
f_{i_1...i_k}
(X_{i_1}^{(1 )},...,X_{i_k}^{(k)}) || \ge t), $

\noindent where again, the last line follows by the decoupling result
for  U-statistics of orders $2,...,k-1$ of the inductive hypothesis.
Since the statement ``not all j's equal'' means that
there are less than $k$ $j$'s equal, the variables whose j's 
are equal can be
decoupled using (conditionally on the other variables)
the decoupling inequalities for U-statistics of order $2,...,k-1$.

Next we give the proof of the lower bound.

\bigskip

\centerline{\bf 5. Proof of the Lower Bound }

\bigskip

In order to show the lower bound we require the following result.

\bigskip

\noindent{\bf Lemma 3.}  Let $1\le l \le k$.  Then there is a constant
$C_k$\ such that
$$ P(||\sum_{1\le i_1\ne i_2\ne ... \ne i_k \le n} f_{i_1...i_k}
   (X_{i_1}^{(1)},X_{i_2}^{(1)},...,X_{i_k}^{(1)}) || \ge  t)\qquad\qquad$$
$$ \qquad\qquad \ge C_k^{-1}
   P(|| \sum_{1\le i_1\ne i_2\ne ... \ne i_k \le n}
   \sum_{1\le j_1,...,j_k\le l} f_{i_1...i_k}
   (X_{i_1}^{(j_1)},X_{i_2}^{(j_2)},...,X_{i_k}^{(j_k)}) || \ge  C_k t)
.$$

\bigskip
\noindent{\bf Proof:}
Let $\{\delta_r\}, r=1,...,l$, 
be a sequence of random variables for which
$P(\delta_r=1) = {1\over l}$ and $  P(\delta_r =0) = 1- {1\over l}$, and
$\sum_{r=1}^l \delta_r = 1$.
Set $\epsilon_r = \delta_r -{1\over l}$ for $ r=1,...,l$.
Then,  it is easy to see that 
there exists $\sigma_l>0$\ depending only upon $l$\
such that for any real number $x_0$ and any
sequence of real constants
$\{a_i\}$
$$|| x_0 +\sum_{r=1}^l a_r \epsilon_r  ||_4 \le 
||x_0+ \sigma_l^{-1} \sum_{r=1}^l a_r \epsilon_r  ||_2.
\leqno (8)$$ 
One can also use the
results of Section 6.9 of Kwapien and Woyczynski (1992) (Pg. 180, 181) 
to assert this since
the $\epsilon $'s satisfy the conditions 1. through 3. stated there.

Let $\{(\delta_{i1},...,\delta_{il}), i=1,...,n\} $ be $n$ independent copies
of 
$(\delta_{1},...,\delta_{l})$. As before, we define 
$$\epsilon_{ij} =
\delta_{ij} - {1 \over l}.\leqno (9)$$
Since the vectors ${\cal E}_i=(\epsilon_{i1},...,\epsilon_{il})$ 
are independent, by an argument given in Kwapien and Szulga (1991), for
$i=1,...,n$,  
for all constants $x_0, a_{ij}$ in $R$,
$$ ||x_0+ \sum_{i=1}^n \sum_{r=1}^l 
 a_{ir}\epsilon_{ir}||_4 \le 
||x_0 + \sigma_l^{-1} \sum_{i=1}^n \sum_{r=1}^l
 a_{ir} \epsilon_{ir} ||_2 \le
\sigma_l^{-1} ||x_0 + \sum_{i=1}^n \sum_{r=1}^l
 a_{ir} \epsilon_{ir} ||_2,\leqno (10)$$
and recentering, we obtain
$$ ||x_0+ \sum_{i=1}^n \sum_{r=1}^l
 a_{ir}\delta_{ir}||_4 \le 
\sigma_l^{-1} ||x_0 + \sum_{i=1}^n \sum_{r=1}^l
 a_{ir} \delta_{ir} ||_2.\leqno (11)$$

Next we use the sequence ${\cal E}_i, i=1,...,n$ in 
defining the analogue of the Z's  
used in our proof of the upper bound.  

For each $i$, let $Z_i= X_i^{(j)}$ if $\delta_{ij} =1$.
Then, $\{Z_i,~i=1,...,n\}$ has the same joint distribution as
$\{X_i^{(1)},~i=1,...,n\}$ and

$$f_{i_1...i_k}(Z_{i_1},..,Z_{i_k}) =
\sum_{1\le j_1, j_2,...,j_k\le l} \delta_{i_1j_1}...\delta_{i_kj_k}
f_{i_1...i_k}(X_{i_1}^{(j_1)},..., X_{i_k}^{(j_k)}).$$

The fact that $E\delta_{i_rj_r} = {1\over l}$ for all $i_r,j_r$ gives, 
$$E(f_{i_1...i_k}(Z_{i_1},...,Z_{i_k})| {\cal G}_l) =
({1\over l})^k\sum_{1\le j_1,...,j_k \le l} 
f_{i_1...i_k}(X_{i_1}^{(j_1)},...,X_{i_k}^{(j_k)}),$$
where ${\cal G}_l = \sigma ((X_i^{(1)}, ..., X_{i}^{(l)} ), ~i=1,...,n).$

\noindent Let

\noindent 
$U_n
=\sum_{1\le i_1\ne i_2 \ne ... \ne i_k \le n} f_{i_1 i_2...i_k} (Z_{i_1},...,
Z_{i_k})$ 
$$
=\sum_{1\le i_1 \ne i_2 \ne ...\ne i_k \le n}\sum_{1\le j_1,...,j_k \le l} 
\delta_{i_1j_1}...\delta_{i_kj_k}
f_{i_1...i_k}(
X_{i_1}^{(j_1)},\cdot \cdot \cdot ,X_{i_k}^{(j_k)}).$$

Let ${\cal D}_i = (\delta_{i1},....,\delta_{il})$.
Since the ${\cal D }$'s are independent of the
$X$'s, if we let 

\noindent
$g_{i_1...i_k}({\bf \cal D}_{i_1},...,{\cal D}_{i_k}) $
$$
=\sum_{1\le j_1,...,j_k \le l} 
\delta_{{i_1}{j_1}}
...\delta_{{i_k}{j_k}} f_{i_1...i_k}(
X_{i_1}^{(j_1)},\cdot \cdot \cdot ,X_{i_k}^{(j_k)}).$$
then, since 
$$f_{i_1...i_k}(X_{i_1},...,X_{i_k})=
f_{i_{(\pi (1))}...i_{(\pi (k))}}(X_{i_{\pi (1)}},...,X_{i_{\pi (k)}}),$$
\noindent we have that,
$$g_{i_1...i_k}({ \cal D}_{i_1},...,{\cal D}_{i_k})=
g_{i_{(\pi (1))}...i_{(\pi (k))}}({\cal D}_{i_{\pi(1)}},...,
{\cal D}_{i_{\pi (k)}}).
$$
Therefore, the two sided decoupling inequality in de la Pe\~na (1992) 
can be applied and, for every convex increasing function $\Phi$,
every ${\cal G}_l$-measurable function $T$,
and $k $ independent copies ${\cal D}_i^{(r)}, ~r=1,...,k$ of 
${\cal D}_i$ there exists numerical constants $A_k$, $B_k$ so that 
$$E (\Phi (A_k||T+\sum_{1\le i_1 \ne i_2 \ne ...\ne i_k \le n}
g_{i_1...i_k}({\cal D}_{i_1},...,{\cal D}_{i_k})||)|{\cal G}_l)$$
$$\le E (\Phi (||T+\sum_{1\le i_1 \ne i_2 \ne ...\ne i_k \le n}
g_{i_1...i_k}({\cal D}_{i_1}^{(1)},...,{\cal D}_{i_k}^{(k)})||)|{\cal G}_l)$$
$$\le E(\Phi (B_k||T+\sum_{1\le i_1 \ne i_2 \ne ...\ne i_k \le n}
g_{i_1...i_k}({\cal D}_{i_1},...,{\cal D}_{i_k})||)|{\cal G}_l).$$

\noindent This result with
(11) shows that conditionally on ${\cal G}_l$
$$||U_n-T_n||_4 \le \sigma_l^{-k} {B_k \over A_k} 
||U_n-T_n||_2, \leqno (12)$$
where
$$T_{n} = E(U_n|{\cal G}_l) =
{({1\over l} )}^k\sum_{1\le i_1 \ne i_2 \ne ... \ne i_k \le n}
\sum_{1\le j_1,...,j_k\le l} f_{i_1...i_k}
(X_{i_1}^{(j_1 )},X_{i_2}^{(j_2)},...,X_{i_k}^{(j_k)}).
$$
(See also the proofs of Lemma 2 and Lemma 6.5.1 of Kwapien and
Woyczynski (1992)).

Thus we have that,
$$P(||U_n || \ge  ||T_n|| |{\cal G}_l) \ge c_k^{-1}.\leqno (13)$$
This follows from the
use of (12) and Proposition 1 with $a=T_n$ and $Y= U_n-T_n$.
We also use the fact that for any
random variable $\xi $ and positive 
constant $c$, $|| \xi ||_4 \le c ||\xi ||_2 $ implies that
$||\xi ||_2 \le c^2 ||\xi ||_1$ 
(See also the proof of Lemma 2
for the approach to transfer the problem from one on
Banach space valued random variables to one on real valued).

Integrating (13) over the set $\{ ||T_n || \ge  t\}$ we get

\noindent$P(||\sum_{1\le i_1 \ne i_2 \ne ... \ne i_k \le n} f_{i_1...i_k}
(X_{i_1}^{(1 )},X_{i_2}^{(1)},...,X_{i_k}^{(1)}) || \ge  t) 
$
$$=P(||\sum_{1\le i_1 \ne i_2 \ne ... \ne i_k \le n} f_{i_1...i_k}
(Z_{i_1},Z_{i_2},...,Z_{i_k}) || \ge  t) 
$$
$$ \ge c_k^{-1}
P( C_k||\sum_{1\le i_1 \ne i_2 \ne ... \ne i_k \le n}
\sum_{1\le j_1,...,j_k\le l} f_{i_1...i_k}
(X_{i_1}^{(j_1 )},X_{i_2}^{(j_2)},...,X_{i_k}^{(j_k)}) || \ge  t),
$$
and Lemma 3 is proved.

The end of the proof of the lower bound follows by using induction and
the iterative procedure introduced to 
obtain the proof of the lower bound multivariate decoupling inequality in
de la Pe\~na (1992). We give a different expression of
the same proof, motivated by ideas from de la Pe\~na, Montgomery-Smith
and Szulga (1992).  We will use $S_k$\
to denote the set of permutations of $\{1,...,k\}$.

The Mazur-Orlicz formula tells us that for any $1 \le j_1,...,j_k \le k$
that 
$$ \sum_{0 \le \delta_1,\dots,\delta_k \le 1}
   (-1)^{k-\delta_1-\dots-\delta_k}
   \delta_{j_1} \dots \delta_{j_k} $$ 
is $0$\ unless $j_1,...,j_k$\ is a permutation of $1,...,k$, in which
case it is $1$. Hence
$$ \eqalignno{
   \sum_{\pi \in S_k}
   f_{i_1\dots i_k}(&X_{i_1}^{(\pi(1))},\dots,X_{i_k}^{(\pi(k))}) \cr
   &= 
   \sum_{0 \le \delta_1,\dots,\delta_k \le 1}
   (-1)^{k-\delta_1-\dots-\delta_k}
   \sum_{1 \le j_1,\dots,j_k \le k}
   \delta_{j_1} \dots \delta_{j_k}
   f_{i_1\dots i_k}(X_{i_1}^{(j_1)},\dots,X_{i_k}^{(j_k)}) . } $$
By the symmetry properties on $f$, 
$$ \eqalignno{
   \sum_{1 \le i_1\ne \dots\ne i_k \le n} &
   f_{i_1\dots i_k}(X_{i_1}^{(1)},\dots,X_{i_k}^{(k)}) \cr
   &= 
   {1\over k!} \sum_{1 \le i_1\ne \dots\ne i_k \le n}
   \sum_{0 \le \delta_1,\dots,\delta_k \le 1}
   (-1)^{k-\delta_1-\dots-\delta_k}
   \sum_{1 \le j_1,\dots,j_k \le k}
   \delta_{j_1} \dots \delta_{j_k}
   f_{i_1\dots i_k}(X_{i_1}^{(j_1)},\dots,X_{i_k}^{(j_k)}) . } $$
Therefore,
$$ \eqalignno{
   \Pr(\| &\sum_{1 \le i_1\ne \dots\ne i_k \le n}
   f_{i_1\dots i_k}(X_{i_1}^{(1)},\dots,X_{i_k}^{(k)}) \| \ge  t) \cr
   &\le
   \sum_{0 \le \delta_1,\dots,\delta_k \le 1}
   \Pr(\| \sum_{1 \le i_1\ne \dots\ne i_k \le n}
   \sum_{1 \le j_1,\dots,j_k \le k}
   \delta_{j_1} \dots \delta_{j_k}
   f_{i_1\dots i_k}(X_{i_1}^{(j_1)},\dots,X_{i_k}^{(j_k)}) \| \ge  
   k!t/2^k) \cr
   &=
   \sum_{l=1}^k \left({k\atop l} \right)
   \Pr(\| \sum_{1 \le i_1\ne \dots\ne i_k \le n}
   \sum_{1 \le j_1,\dots,j_k \le l}
   f_{i_1\dots i_k}(X_{i_1}^{(j_1)},\dots,X_{i_k}^{(j_k)}) \| \ge  
   k!t/2^k) ,\cr } $$
and this combined with Lemma~3 is sufficient to show the result.

\vfill\eject

\centerline {\bf7. References}

\item {\bf 1.} de la Pe\~na, V. H. (1992). Decoupling and Khintchine's 
inequalities for U-statistics. {\it Ann. Probab.} {\bf 20} 4, 1877-1892.

\item {\bf 2.} de la Pe\~na, V. H., Montgomery-Smith 
and Szulga, J. (1992). Contraction and decoupling 
inequalities for multilinear forms and U-statistics. Preprint.

\item {\bf 3.} de la Pe\~na, V. H. and Montgomery-Smith (1993).  
Bounds on the tail probability of U-statistics and quadratic forms.
Preprint.

\item{\bf 4.} Gin\'e, E. and Zinn, J. (1992). 
A remark on convergence in distribution of U-statistics.
Preprint.

\item {\bf 5.} Kwapien, S. Szulga, J. (1991). Hypercontraction methods
in moment inequalities for series of independent random variables in
normed spaces.
{\it Ann. Probab.} {\bf 19} (1), 369-379.

\item {\bf 6.} Kwapien, S. and Woyczynski, W. (1992). Random Series and 
Stochastic Integrals: Single and Multiple. Birkhauser.

\item{\bf 7.} McConnell, T. and Taqqu, M. (1986).
Decoupling inequalities for multilinear forms in independent symmetric random
variables. {\it Ann. Probab.} {\bf 14}(3), 943-954.

\bye
\vfill\eject

The Mazur-Orlicz formula may be rewritten thus:
$$ \eqalignno{
   \sum_{\pi \in S_k}
   f_{i_1\dots i_k}(&X_{i_1}^{(\pi(1))},\dots,X_{i_k}^{(\pi(k))}) \cr
   &= 
   \sum_{0 \le \delta_1,\dots,\delta_k \le 1}
   (-1)^{k-\delta_1-\dots-\delta_k}
   \sum_{1 \le j_1,\dots,j_k \le k}
   \delta_{j_1} \dots \delta_{j_k}
   f_{i_1\dots i_k}(X_{i_1}^{(j_1)},\dots,X_{i_k}^{(j_k)}) . } $$
By the symmetry properties on $f$, 
$$ \eqalignno{
   \sum_{1 \le i_1\ne \dots\ne i_k \le n} &
   f_{i_1\dots i_k}(X_{i_1}^{(1)},\dots,X_{i_k}^{(k)}) \cr
   &= 
   {1\over k!} \sum_{1 \le i_1\ne \dots\ne i_k \le n}
   \sum_{0 \le \delta_1,\dots,\delta_k \le 1}
   (-1)^{k-\delta_1-\dots-\delta_k}
   \sum_{1 \le j_1,\dots,j_k \le k}
   \delta_{j_1} \dots \delta_{j_k}
   f_{i_1\dots i_k}(X_{i_1}^{(j_1)},\dots,X_{i_k}^{(j_k)}) . } $$
Therefore,
$$ \eqalignno{
   \Pr(\| &\sum_{1 \le i_1\ne \dots\ne i_k \le n}
   f_{i_1\dots i_k}(X_{i_1}^{(1)},\dots,X_{i_k}^{(k)}) \| \ge  t) \cr
   &\le
   \sum_{0 \le \delta_1,\dots,\delta_k \le 1}
   \Pr(\| \sum_{1 \le i_1\ne \dots\ne i_k \le n}
   \sum_{1 \le j_1,\dots,j_k \le k}
   \delta_{j_1} \dots \delta_{j_k}
   f_{i_1\dots i_k}(X_{i_1}^{(j_1)},\dots,X_{i_k}^{(j_k)}) \| \ge  k!t/2^k) \cr
   &=
   \sum_{l=1}^k \left({k\atop l} \right)
   \Pr(\| \sum_{1 \le i_1\ne \dots\ne i_k \le n}
   \sum_{1 \le j_1,\dots,j_k \le l}
   f_{i_1\dots i_k}(X_{i_1}^{(j_1)},\dots,X_{i_k}^{(j_k)}) \| \ge  k!t/2^k) .\cr
    } $$

\bye